\documentclass{amsart}

\usepackage{mathrsfs} 
\usepackage{amsmath}
\usepackage{algorithm}
\usepackage[noend]{algpseudocode}
\usepackage{graphicx}
\usepackage{color}

\makeatletter
\def\BState{\State\hskip-\ALG@thistlm}
\makeatother

\begin{document}

	\title{Optimizing the Production Cost of Minting with Mixed Integer Programming}
	\author{Carlos A. Alfaro}
	\address{Banco de M\'exico, Mexico City 11500, Mexico}
	\email{carlos.alfaro@banxico.org.mx}

	\author{Ra\'ul Mart\'inez-Noriega}
	\address{Banco de M\'exico, Mexico City 11500, Mexico}
	\email{raulm.noriega@banxico.org.mx}

	\author{C\'esar Guadarrama}
	\address{Banco de M\'exico, Mexico City 11500, Mexico}
	\email{cesar.uribe@banxico.org.mx}

	\author{Adolfo S\'anchez-Flores}
	\address{Banco de M\'exico, Mexico City 11500, Mexico}
	\email{asanchef@banxico.org.mx}

	\author{Jorge A. Aguilera}
	\address{Banco de M\'exico, Mexico City 11500, Mexico}
	\email{jaguilera@banxico.org.mx}

\maketitle

\markleft{CARLOS A. ALFARO ET AL.}

\begin{abstract}
	For central banks, managing the minting is one of the most important task since a shortage yields negative economic and social impacts, and the budget committed for minting is one of the largest within the central banks. 
	Hence, the central bank requires to find the mixture of coins to be produced that satisfies the demand, inventory and production constraints while minimizing the cost. 
	We propose a mixed-integer programming model that minimize the cost of minting by reducing the number of extra-shifts required while fulfilling the constraints.  
	We also perform a simulation with data of a central bank which shows that the model reduces in 24\% the cost of extra-shifts used during 21 quarters, compared with the spreadsheet based approach used currently at the operation.
\end{abstract}
		



\section{Introduction}
One of the most important roles for any central bank (CB) is to assure a sufficient and uninterrupted supply of cash to match the demand generated by the country's economic activity.  
A shortage of cash would impose constraints on the economy by reducing the possibility of doing basic transactions, besides annoyance among the population and services. 
Despite the availability of alternative payment methods, such as credit, debit, and prepaid cards, on-line banking, and most recently mobile banking, cash is still the favorite one in development countries.
For instance, in Mexico \cite{Mazzotta} \%90 of the consumer transactions are still performed in cash.
Therefore, managing of cash production becomes a key aspect for any nation.

The processes involved in satisfying the demand of cash are forecasting of the demand, definition of the safety stock levels, fabrication process and, distribution. 
In this study, we will focus in minting of circulating coins, which is the process where the coins are produced; more specifically, we are interested in minimizing the minting cost.

Minting process is usually the most costly due to the high cost of raw materials. 
The country's minting policies define the coin characteristics and its cost, these might vary between countries.
A minting order consists of the number of coins to be produced for each denomination in a specific quarter.
In the Mexican case, central bank of Mexico (CBM) is in charge of submitting the minting order to the Mexican mint (MM), who will produce the coins.
For this minting order, CBM pays to MM the metal and workforce required to manufacture it.

The manufacture process is composed of three processes: \emph{blanking}, \emph{annealing} and \emph{striking}. 
Blanking is the process where the coin's pieces are cutting out of metal sheets. 
Then, the annealing process starts, in which those pieces are heated to make uniform the hardness of the coin metals. 
Finally, striking is the process where the motives are stamped in both faces of the coin.
The cost function for each minting process consists a piecewise linear function.

Therefore, the central bank requires to find the mixture of coins (minting order) to be produced at each quarter that minimize the cost while satisfying the demand, inventory and production constraints.
In practice, this problem is usually solved manually in a stepwise procedure.

The main contribution of this article is to present a mathematical model formulation to reduce the cost of minting  by minimizing the number of extra-shifts required during the yearly production, while assuming real world constraints. 
Additionally, we develop a metaheuristic solution method based on greedy heuristics and a relaxation of the MIP model.
A sequence of computational simulations with this solution is performed using historical data of 21 quarters of a central bank.
The presented solution approach and its implementation is explained in detail and tested in extensive numerical simulation.
This simulation was implemented in C++ while the MIP model was solved using the GLPK solver \cite{GLPK}.
By comparing the accumulated costs between the simulation and the historical solution, it is shown that the model reduces to 24\% the cost of extra-shifts used during 21 quarters of the observed data.

The article is organized as follows: in Section~\ref{sec:literature}, an overview of the current literature on this problem is given.
Section~\ref{sec:Definicion} refers to the problem presentation and definition of the minting cost functions.  
The model and its constrains are described in Section~\ref{sec:Optimizacion}. 
In Section~\ref{sec_simulation} the greedy heuristics are introduced to solve the problem, together with the simulation of a real world scenario. 
Finally we conclude this paper in Section~\ref{sec_conclusion}. 

\section{Literature review}\label{sec:literature}
The literature dealing with minting cost optimization is rather scarce. 
However there are plenty papers dealing with different aspects in satisfying the coins demand. 

In \cite{Chumacero_2007} there were reviewed ARIMA and VaR models for forecasting the coin and bank note demand. 
And in \cite{Venkatesh_2014} the authors used clustering and neural networks for predicting the cash demands in automatic teller machines (ATMs). 

Optimizing the stock levels of banknotes has been reported by Baumol in ~\cite{Baumol_1952}, where he proposed a simple static bank note inventory model. 
Ladany~\cite{Ladany_1997} proposed in 1997 a discrete dynamic programming model for the note ordering policy at the Bank of Israel. 
More recently, in 2005, Massoud~\cite{Massoud} provided a dynamic cost minimizing the banknote inventory model applied to a scenario similar to the one of the Central Bank of Canada.

Optimizations related with the banknote distribution are reported by Castro,~\cite{Castro_2009}, in 2009 and Osorio \& Toro,~\cite{OT_2012}, in 2012.
Both proposals are related to the distribution of notes to ATMs using stochastic programming models. 
Recently, in 2015, Zhu et. al.~\cite{Zhu_2015} researched the logistic problem faced by certain central bank to supply cash to its regional branches under assumptions of security concerns.

The actual process of manufacturing coins or minting has been researched in~\cite{Yen_1981} and~\cite{Leitao_97}. 
Another interesting point is studied in \cite{Bouhdaoui_2011} and \cite{VanHove_1996}, where it is shown that the spacing between denominations can increase the production cost incurred by the CB.

A proposal describing the optimization of minting is reported by Ladany~\cite{Ladany_1981}. 
He proposed a model to control the minting of commemorative coins to obtain the maximum profit when the coins are sold, e.g. to collectors. 

\section{Problem description}\label{sec:Definicion}
Based on the demand, the CB emits a \emph{minting plan}  $\mathbf{F}$ that describes the mixture of coins to be produced by the Mint during the planning horizon $\mathcal{T}$. 
The minting plan $\mathbf{F}=(\mathbf{f}_t)_{t\in \mathcal{T}}$ consists of quarterly \emph{minting orders} $\mathbf{f}_t=(f_t^d)_{d\in\mathcal{D}}$, which details the number of coins $f_t^d$ per denomination $d$ that the Mint should manufacture at quarter $t$. 
Then, the Mint computes the requirements to produce those quantities and set a cost that is covered by the CB. 

The CB pays to the Mint a fixed rent for its services. 
The rent cost includes the basic level production capacity for each process, \emph{i.e.}, blanking, annealing and striking. 
An extra capacity increases the quarterly production to different levels depending on the configuration. 
These production levels can be hired, individually for each process, depending on the demand requirements.
However, any extra capacity has an additional cost to the CB, this is because any extra capacity requires to use third shifts. 
Thus, depending on the mixture of coins, the Mint may require extra-shifts to fulfill the order, incurring in extra costs to the CB.

The precise production level to be hired for each process depends on the mixture, i.e. the number of coins to manufacture of each denomination. 
This is due to coins of different denomination vary in size and weight, and bigger coins require more resources of particular processes. 
Moreover, there are mono-metallic (single piece) and bi-metallic coins (constructed by a core of certain metal and ringed by another metal) which demand different amount of the processes.

Then, our objective is given the initial coin inventory levels, coin demand and mint production restrictions, to find a suitable set of minting orders that satisfies the quarterly demand that minimize the use of the extra-shifts.

\subsection{Assumptions}
Issuing an optimal demand for coins is not directly pursued since we deal with the minimization of the production cost of current coins given any demand and forecasting error.
In our study, we do not contemplate special commemorative coins.
Also, we are not interested in optimizing the manufacturing process, we rather optimize on the mixture of coins satisfying manufacturing processes fixed by the mint.
This is because although the minting is managed by the CB, the minting is conveyed by the mint which is an independent institution to the CB.
In our proposal, the inventory levels are already defined and fixed by the CB's policies. 
Also, we assume that the spacing between denominations was previously fixed. 

\subsection{Minting processes and their cost}

In the following, we introduce the cost functions associated to each process, blanking, annealing, and striking. These functions were defined based on the contracts that the CB holds with the Mint.

\subsubsection{Blanking}
Blanking is the first process of the production line and consists in obtaining the coin pieces out of metal sheets. 

The production level for this process is reported by the Mint as the number of coins of certain denomination that can be processed in one day.  
Given a mintage order $\mathbf{f}_t$ we can measure the use of blanking by transforming the number of coins of each denomination into working days and then add them up. 
Let $D(\mathbf{f}_t)$ denote the transformation that convert a minting order $\mathbf{f}_t$ into the number of working days required to be produced.
This is dot product between $\mathbf{f}_t$ and a constant vector.

Equation~(\ref{costoD}) is the cost function for blanking. 
The first row refers to the base capacity and its cost is zero because it is already included in the rent. 
From the second row onwards are the costs of different configurations for the extended capacity.
\begin{equation}\label{costoD}
	\mathscr{C}(\mathbf{f}_t) = 
	\begin{cases}
   				0, & \text{if } D(\mathbf{f}_t) \leq x_0,\\
   				C_1, &	\text{if } x_0 < D(\mathbf{f}_t) \leq x_1,\\
   				\vdots\\
   				C_{nc} & \text{if } x_{nc-1} < D(\mathbf{f}_t) \leq x_{nc},\\
  	\end{cases}
\end{equation}
where $C_1,\dots,C_{nc}$ are the costs associated to certain production level, and $x_i$ denote the number of days required.

\subsubsection{Annealing}
Annealing is the process where the coins receive a heat treatment to uniform the hardness of the pieces based on copper alloys. Only the bi-metallic coins uses this alloy.

The production capacity of this resource depends on the weight of copper alloy pieces to furnace in a quarter, and therefore this depends on the number of coins per denomination defined in the minting order. 
Assuming that $w^d$ is the weight of the copper alloy corresponding to the denomination $d$, then the total weight $W(\mathbf{f}_t)$ of a minting order $\mathbf{f}_t$ is
$$
	W(\mathbf{f}_t) = \sum_d{w^df^d_t}.
$$

There is only one extra-shift level after the basic capacity for annealing. 
Thus, the cost function for this process is
\begin{equation}\label{costoH}
	\mathscr{H}(\mathbf{f}_t) = \left\{
  \begin{array}{l l l}
   0, &\quad \text{if} &\quad  W(\mathbf{f}_t) \leq y_0\;\text{Ton.}\\
   H, &\quad \text{if} &\quad y_0 < W(\mathbf{f}_t) \leq y_1 \;\text{Ton.,}
  \end{array} \right. 
\end{equation}
where $H$ is the cost of annealing more than $y_0$ Tons of coins.
In this case we also consider the base capacity to have a zero cost because it is included in the rent.

\subsubsection{Striking}
Striking is the last process and it is where the coin's faces are stamped. In this case, the production capacity is independent of the denomination and only depends on the number $N$ of minted coins in a quarter, where
$$
	N({\bf f}_t) = \sum_d{f^d_t}.
$$

Equation~(\ref{costoA}) describes the striking cost function. The first row represent the base capacity and the rows what follows, describe different configurations of extended capacity. 
\begin{equation}\label{costoA}
	\mathscr{A}({\bf f}_t) = \left\{
	\begin{array}{l l l}
		0, &\; \text{if} &\; N({\bf f}_t) \leq z_0\\
   		A_1, &\; \text{if} &\; z_0< N({\bf f}_t) \leq z_1\\
   		\vdots\\
   		A_{na},&\; \text{if}  &\;  z_{na-1}< N({\bf f}_t) \leq z_{na}
  \end{array} 
	\right. ,
\end{equation}
where $A_1,\dots,A_{na}$ are the costs of extended capacities.

\section{Mathematical model}\label{sec:Optimizacion}
The objective is to minimize the annual production cost, that is, reducing the cost of the minting orders. 
According to Equations~(\ref{costoD}),~(\ref{costoH}), and (\ref{costoA}), reducing the cost of a minting order is equivalent to minimize the use of extra shifts at each production process.  
In this section we propose a model that minimize the usage of  extra shifts for a planning horizon of $T$ quarters. 
Solving the model implies finding the minting orders that minimize the cost and meet all inventory and production constraints.

We begin by defining the total production cost for a planning horizon of $T$ quarters as the sum of all production costs used during the planning horizon:
\begin{equation}\label{seudoObjective}
	\sum_t{\sum_i{C_i c_t^i }}+\sum_t{Hh_t} + \sum_t{\sum_j{A_j a_t^j }} ,
\end{equation}
where $t=1,2,\dots,T$, $i=1,2,\dots,nc$, $j=1,2,\dots,na$, and $c^i,h,a^j$ are binary decision variables used to select a cost represented by $C_i$, $H$, $A_j$ and described in~(\ref{costoD}), (\ref{costoH}), (\ref{costoA}) respectively.

\begin{table}[h]
\centering
\begin{tabular}{c@{\extracolsep{.5cm}}l}
\hline
 & \bf Indices \\
 \hline
 $i$ & third-shift level of blanking\\
 $j$ & third-shift level of striking\\
 $t$ & quarter index $t\in\{1, \dots, T\}$\\
 $d$ & denomination\\
 & \\
 \hline
  & \bf Parameters \\
 \hline
 $T$ & number of quarters in the planning horizon\\
 $C_i$ & cost of blanking third-shift at level $i$ \\
 $A_j$ & cost of striking third-shift at level $j$ \\
 $H$ & cost of the annealing third-shift\\
 $P_t^d$ & coin demand for denomination $d$ at quarter $t$\\
 $DEM_t^d$ & average coin demand of three months for denomination $d$ at quarter $t$\\
 $IMAX$ & maximum vault capacity\\
 $IMIN^d$ & minimal inventory of denomination $d$\\ 
 $x_i$ & upper bound for the number of working days required for third-shift level $i$\\
& of blanking\\
 $y_1$ & upper bound for the weight required for the third-shift of annealing\\
 $z_j$ & upper bound for the number of minted coins required for third-shift level $j$\\
& of striking\\ 
& \\
 \hline
  & \bf Variables\\
 \hline
 $c_t^i$ & binary variable represents the assigning of third-shift level $i$ of blanking\\
 $h_t$ &binary variable represents the assigning of the third-shift of annealing\\
 $a_t^j$ &binary variable represents the assigning of third-shift level $j$ of striking\\
 $f_t^d$ & integer variables representing the number of minted coins of denomination $d$\\
&  at quarter $t$\\
 $E_t^d$ & integer variables representing the coin inventory of denomination $d$\\
&  at quarter $t$\\
 $K$ & real positive variable used to efficiently use the production levels\\
 \hline\\
\end{tabular}
\caption{Summary of notations}
\end{table}

If we consider~(\ref{seudoObjective}) as the objective function to minimize, the solution will provide us the production levels for each process that satisfies the demand and production's constrains. 
However, those levels can be misused because this solution does not maximize the number of coins to be produced. 
We deal with this issue by adding a variable $K$  in~(\ref{seudoObjective}) to maximize the number of coins, where $K$ is related to the stocks by means of the constraint~(\ref{eq_rProduccion}). 
In practice $K\sim1$, since the existence at the end of the planning horizon can not be much larger than the minimum inventory.
Thus, our objective function to minimize is
\begin{equation} \label{eq:funcObjetivo}
	\sum_t{\sum_i{C_i c_t^i }}+\sum_t{Hh_t} + \sum_t{\sum_j{A_j a_t^j }}   - K.
\end{equation}
Therefore this bi-objective function aims minimizing the extra cost of the production process and maximizing the number of coins to produce, giving priority to minimize the cost. 

Then the objective function,~(\ref{eq:funcObjetivo}), is solved holding the following constrains:
\begin{align}
	&\forall t, & W(\mathbf{f}_t) \leq y_0 + (y_1-y_0)h_t \label{eq_rHorno}\\
	&\forall t, & \sum_d f^d_t \leq z_0 + \sum_j{(z_j-z_{j-1}) a_t^j } \label{eq_rAcunia}\\
	&\forall t, & \sum_j a_t^j \leq 1 \label{eq_rAcunia_1}\\
	&\forall t, & D(\mathbf{f}_t) \leq x_0 + \sum_i{(x_i-x_{i-1}) c_t^i } \label{eq_rCorte}\\
	&\forall t, & \sum_i c_t^i \leq 1 \label{eq_rCorte_1}\\
	&\forall t,d, & E_t^d = E^d_{t-1}+f^d_t-P^d_t \label{eq_rExistencias}\\
	&\forall d, & E_T^d \geq IMIN^d K \label{eq_rProduccion}\\
	&\forall t,	& \sum_d{E^d_t}\leq IMAX \label{eq_rIMAX}\\
	&\forall t,d,	& E^d_t \geq DEM^d_t \label{eq_rIMIN}
\end{align}

The first constrain,~(\ref{eq_rHorno}), means that the total weight of the copper alloys that form a minting order should be lower than the base capacity, $(h_t=0)$, or than the extended capacity, $(h_t=1)$. 
Inequalities~(\ref{eq_rAcunia}) and~(\ref{eq_rCorte}) guarantee that the number of coins to strike or to anneal, respectively, are into a feasible production levels at each quarter. 
Then,~(\ref{eq_rAcunia_1}) and~(\ref{eq_rCorte_1}) secure to evaluate one production level at a time.

The equality~(\ref{eq_rExistencias}) relates the expected stock $E^d_t$ for each denomination $d$ with the forecasting demand  $P^d_t$, the production $f^d_t$, and the initial stock level  $E^d_{t-1}$ at each quarter $t$. 
Constraint (\ref{eq_rIMIN}) forces the inventory $E^d_{t}$ is greater or equal than the average demand of three months $DEM_t^d$, this provides an appropriated inventory level to operate.
Constraint (\ref{eq_rIMAX}) guarantee that the total inventory at each quarter is no larger than the maximum capacity $IMAX$.

As mentioned before,~(\ref{eq:funcObjetivo}) aims a double optimization, the minimum overall production cost and then maximization of the number of coins. 
The latter is achieved with the variable $K$ and the constraint~(\ref{eq_rProduccion}). 
This constraint increases the production uniformly, in terms of the minimal inventory of denomination $d$, $IMIN^d$.

The solution to the model are obtained such that~(\ref{eq:funcObjetivo}) yields the minimum cost. 
Specifically, the solution $\sf{S}=(\sf{O},\sf{F},\sf{A},\sf{C},\sf{H})$ provides the optimum cost $\sf{O}$, the minting orders $\sf{F}=\{f_t^d\}$, the used striking extra-shift levels $\sf{A}=\{a_t^j\}$, the used blanking extra-shift levels $\sf{C}=\{c_t^i\}$ and the used annealing extra-shift levels $\sf{H}=\{h_t\}$.
	
\section{Solution approach and simulation}\label{sec_simulation}
This section describes some extra considerations that are needed to solve the model in a real world scenario.
We compare the solution of the MIP model described in Sec.~\ref{sec:Optimizacion} with 21 historical minting orders performed by the CB in the past.

The historical minting orders are proprietary data of the CB and due to privacy reasons we cannot disclose all the details. 
All the results reported in this section were transformed to keep their privacy.

\subsection{Data description}
\begin{figure}[h!]
	\centering
	\includegraphics[width=12cm]{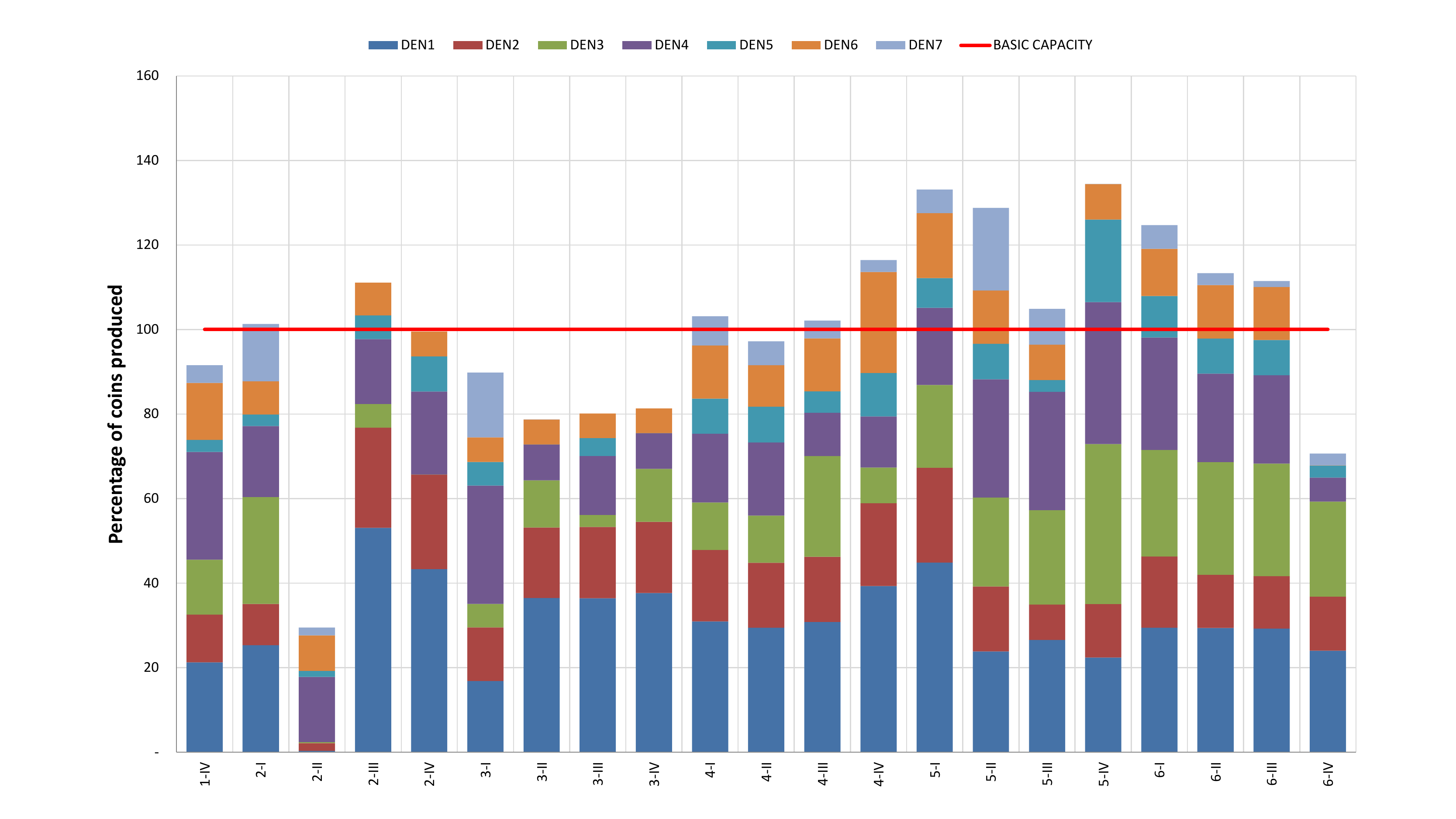}
	\caption{Historical minting orders. The orders consist of the number of pieces produced for 7 denominations. The red line serves as reference to the base capacity of striking.}
	\label{fig:MintageOrdersObserved}
\end{figure}

The data includes observed and forecasting demand, vault and inventory restrictions, and the historical mintage orders during 21 consecutive quarters. 
For example, Fig.~\ref{fig:MintageOrdersObserved} shows the historical minting orders expressed in millions of coins with the $x$-axis representing the quarters in a year-quarter fashion. 
The striking base capacity is shown with a horizontal red line. 
Notice that the quarter 2-II has low production level, the reason was a supply problem of gas which limited the production. 
We will maintain the same scenario, adding the gas supply problem as a restriction to the MIP model in that quarter, to make the proposed model comparable with the historical minting orders.

\subsection{Solution}
To make the problem tractable, we relax model in order that the originally integer variables $f_t^d$ and $E_t^d$ to be real positive.
Then, a solution of this relaxation of the MIP model was obtained using the MIP solver GLPK, which uses Branch \& Cut as primary method to solve the model \cite{GLPK}. 
In general the scenarios were able to solve them in a reasonable computer-time, around 2 seconds each instance.
In addition, the simulation was executed using meta-heuristics that deal with the lack of long-term information and with forecasting errors.
After the solution of the relaxation were obtained, we perform a local-search to obtain a near-optimal solution to the original model.

Now we will introduce a few simple meta-heuristics to improve the iterative optimization of the scenarios. 
These can be reinterpreted as soft constrains added to the model when additional manipulation is required. 

{\bf{Procedure 1.}}
Sometimes the model's solution suggests that the minting order for certain planing horizon $T$ should be produced using least production capacity than the base capacity in any of the three production processes. 
This result is motivated by the combinations of stock's levels and forecasting demand, whose expectation is that there is not need of more production. 
However, an unexpected increase in the future demand could trigger the use of extended capacity in the near future, which would substantially increases the cost. 
For those cases we push the model to fully use the base capacity if and only if this increment does not require of extended capacity in any of the other processes.

We consider this adjustment only for blanking and striking since they are the tighter and most expensive processes. 
This procedure is detailed in Algorithm~\ref{heuristic1} where the following constraints are added to the model:
\begin{equation}\label{Eqn:rule1}
	\text{(striking) }\sum_d f^d_1 = z_0  \text{\hspace{1cm} or \hspace{1cm} (blanking) } D(\mathbf{f}_1) = x_0.
\end{equation}

\begin{algorithm}[h]
\caption{}\label{heuristic1}
\begin{algorithmic}[1]
\BState \textbf{Input}: A solution $\sf{S}=(\sf{O},\sf{F},\sf{A},\sf{C},\sf{H})$ of the model $\sf M$ proposed in Sec.~\ref{sec:Optimizacion}.
\If {$\sum_d {\sf f}^d_1 < z_0$}
\State Let $\sf{S}'=(\sf{O}',\sf{F}',\sf{A}',\sf{C}',\sf{H}')$ be the solution of the model $\sf M'$ obtained by adding the striking constraint of (\ref{Eqn:rule1}) to $\sf M$.
\If {$\sf O=O'$}
\State Take ${\sf M\leftarrow M'}$ and $\sf S\leftarrow S'$.
\EndIf
\EndIf
\If {$D({\sf f}_1) < x_0$}
\State Let $\sf{S}'=(\sf{O}',\sf{F}',\sf{A}',\sf{C}',\sf{H}')$ be the solution of the model $\sf M'$ obtained by adding the blanking constraint of (\ref{Eqn:rule1}) to $\sf M$.
\If {$\sf O=O'$}
\State Take ${\sf M\leftarrow M'}$ and $\sf S\leftarrow S'$.
\EndIf
\EndIf
\State\Return $\sf S$.
\end{algorithmic}
\end{algorithm}

{\bf{Procedure 2.}}
Another interesting situation is when the solution requires the extended capacity at the current quarter. In those cases it is desirable to postpone the use of extended capacity to the next quarter if and only if the solution is still feasible. The motivation for this rule is that the extended capacity could be triggered by a forecasting error on the demand. Therefore, by postponing the use of extended capacity we avoid the use of resources that are not completely needed at the moment, and save their associated costs. On the other hand if the extended capacity is still needed, a major utilization over the extended capacity could be done obtaining a reduction on the unitary cost. 
This procedure is detailed in Algorithm~\ref{heuristic2} where the following constraints are added to the model:

\begin{equation}\label{Eqn:rule2}
	\text{(striking) }\sum_j a_1^j = 0, \hspace{0.5cm} \text{(blanking) }\sum_i c_1^i = 0, \hspace{0.5cm}\text{ or \hspace{0.5cm}(annealing) } h_1 = 0.
\end{equation}

\begin{algorithm}[h]
\caption{}\label{heuristic2}
\begin{algorithmic}[1]
\BState \textbf{Input}: A solution $\sf{S}=(\sf{O},\sf{F},\sf{A},\sf{C},\sf{H})$ of the model $\sf M$ proposed in Sec.~\ref{sec:Optimizacion}.
\If {$\sum_d {\sf f}^d_1 > z_0$}
\State Let $\sf{S}'=(\sf{O}',\sf{F}',\sf{A}',\sf{C}',\sf{H}')$ be the solution of the model $\sf M'$ obtained by adding the striking constraint of (\ref{Eqn:rule2}) to $\sf M$.
\If {$\sf S'\neq \emptyset$}
\State Take ${\sf M\leftarrow M'}$ and $\sf S\leftarrow S'$.
\EndIf
\EndIf
\If {$D({\sf f}_1) > x_0$}
\State Let $\sf{S}'=(\sf{O}',\sf{F}',\sf{A}',\sf{C}',\sf{H}')$ be the solution of the model $\sf M'$ obtained by adding the blanking constraint of (\ref{Eqn:rule2}) to $\sf M$.
\If {$\sf S'\neq\emptyset$}
\State Take ${\sf M\leftarrow M'}$ and $\sf S\leftarrow S'$.
\EndIf
\EndIf
\If {$W({\sf f}_1) > y_0$}
\State Let $\sf{S}'=(\sf{O}',\sf{F}',\sf{A}',\sf{C}',\sf{H}')$ be the solution of the model $\sf M'$ obtained by adding the annealing constraint of (\ref{Eqn:rule2}) to $\sf M$.
\If {$\sf S'\neq\emptyset$}
\State Take ${\sf M\leftarrow M'}$ and $\sf S\leftarrow S'$.
\EndIf
\EndIf
\State\Return $\sf S$.
\end{algorithmic}
\end{algorithm}

\subsection{Simulation}
The annual minting plan is defined for the first time by the CB at the third quarter of the year and, due to CB's policies, the planning horizon of the minting plan starts with the four quarter of the year and must contain the four quarter of the next year. 
After, the information such as observed demand, forecasted demand and inventory levels are updated each quarter, until the third quarter is reached. 

Therefore, the simulation was carried on step by step updating information and solving the optimization at each quarter. 
Initially, the minting plan is defined from four quarter, and the optimization is solved with a planning horizon $T=5$, i.e. solving the minting order of the last quarter in the current year and all minting orders of the next year. Then, the next quarter the information is updated and the optimization is solved with $T=4$, and so on. This process of solving-and-updating is repeated until the scenario with $T=2$, after that, a new minting plan should be determined.

This iterative fashion of solving the proposed model may dismiss the opportunity costs in comparison with a long-term optimization because when $T$ is small, the model cannot consider the future demand and the base capacity can be misused or we can lost the chance to postpone the use of extended capacity.
However, this is the way in which the current CB's policies define the process.
Here is were the heuristics introduced in Subsection 5.2 become worthy, since they help us to deal with these inconveniences of the coin policy. 


\subsection{Results}
We compare the minting orders obtaining with the proposed model's simulation and the historical minting orders of 21 consecutive quarters made by the CB. 

\begin{figure}[h!]
	\centering
	\includegraphics[width=12cm]{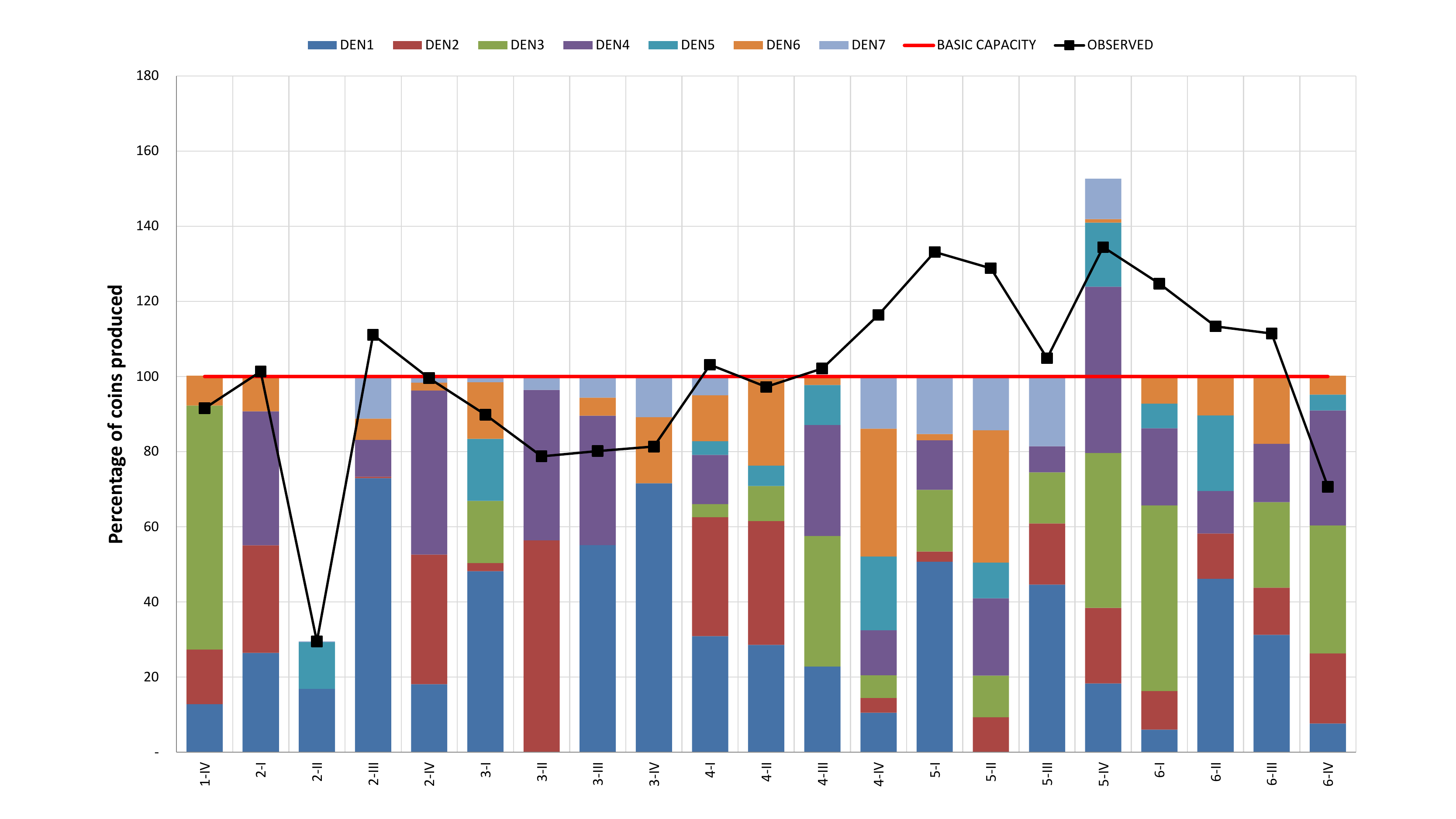}
	\caption{Simulation mintage orders. The orders consist of the number of pieces produced of the 7 denominations in a determined quarter. The red line serves as reference to the maximum basic capacity of the striking. The black line shows the total number of pieces of the Historic mintage orders.}
	\label{fig:MintageOrdersSimulation}
\end{figure}

Fig.~\ref{fig:MintageOrdersSimulation} shows the quarterly minting orders obtained with the proposed model. Recall that at each quarter or scenario, the solutions give us $T$ minting orders corresponding to the planning horizon. However Fig.~\ref{fig:MintageOrdersSimulation} only shows the first order corresponding to certain quarter, the remaining orders are not shown because they will be computed again at the  following quarters with updated information. The black line is the number of coins produced in the historical minting orders. We also add a red line as reference to the striking base capacity which is the most costly process. 

We observe that, with the proposed model, the striking base capacity was completely used in all quarters excepting the 2-II in which there was a supply problem of gas at the Mint and the production was reduced. Using full base production capacity at each quarter, which has no extra cost for the CB except for the raw materials used to mint, allow us to reduce the times that extended capacity was used and thus reducing the cost. During the 21 quarters, the proposed model used only once the extended capacity for striking, at 5-II quarter, compared with the historical orders where the extended capacity was contracted 12 times as side effect of misusing the base capacity in 7 quarters. We showed only striking production capacity because it is the process that uses the most extra-shifts and therefore the most costly.

\begin{figure}[h!]
	\centering
	\includegraphics[width=12cm]{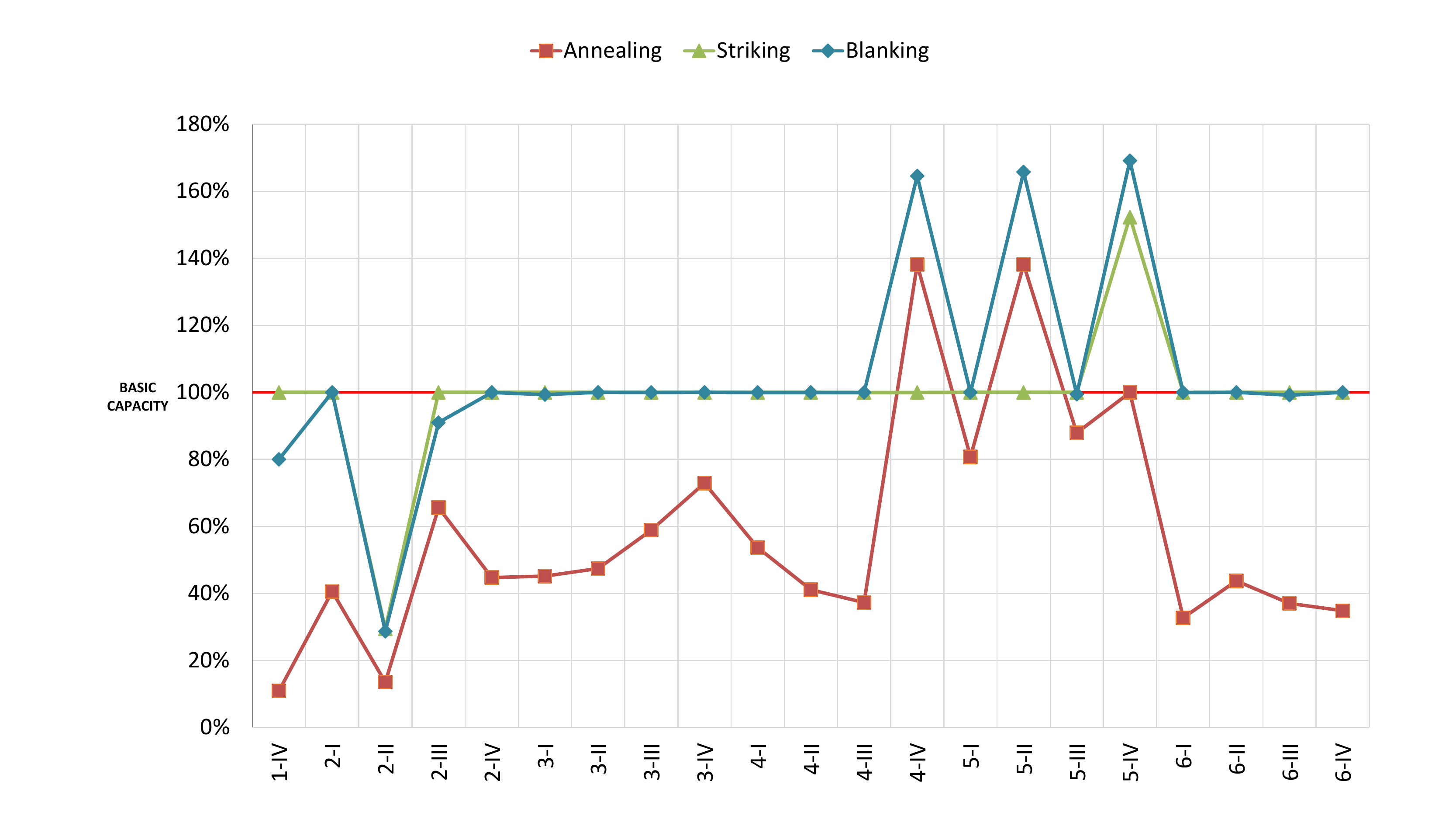} \\
	(a) Simulation\\
	\includegraphics[width=12cm]{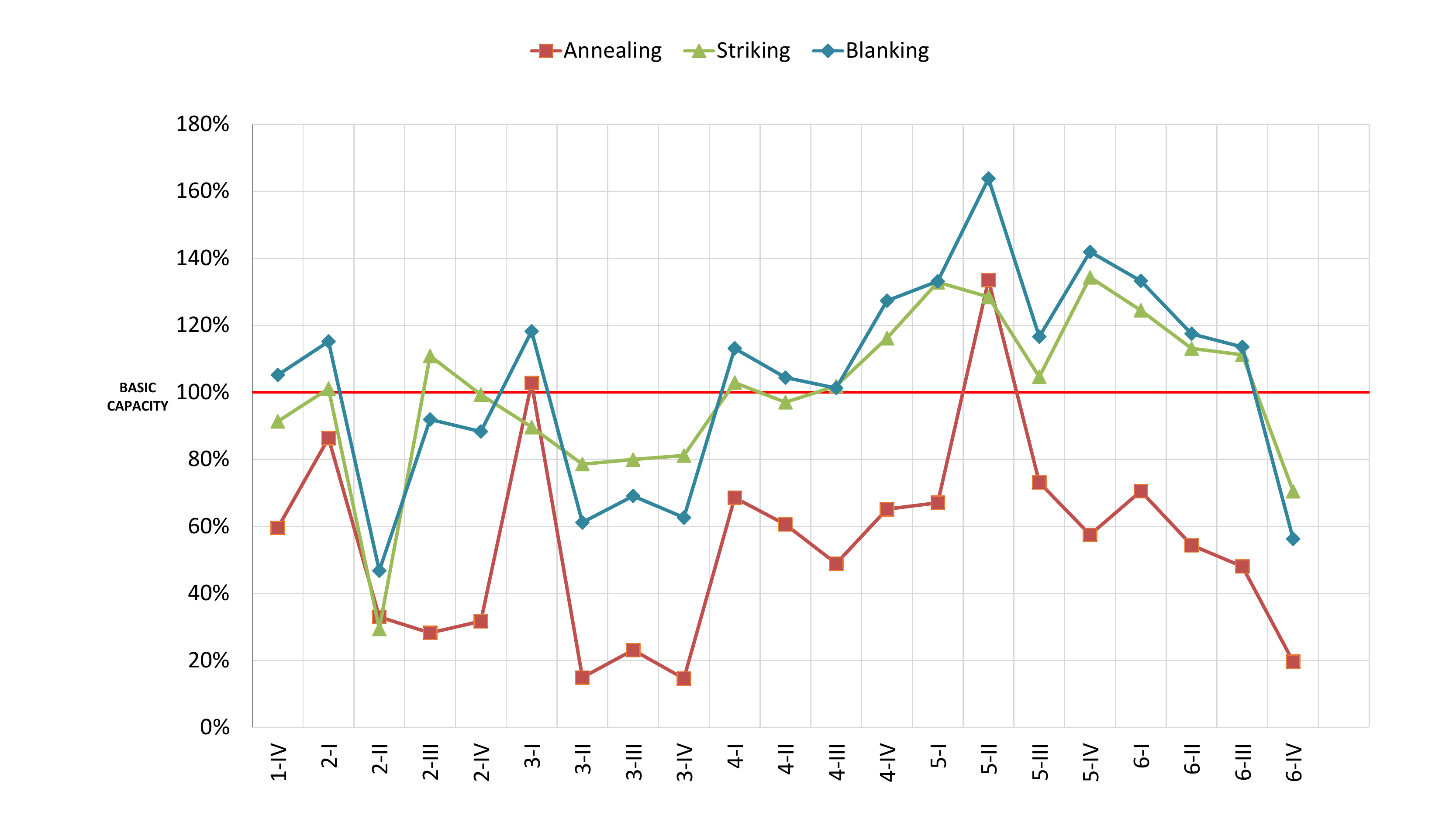}\\
	(b) Observed\\
\caption{Percentage of the utilization of the processes.}
\label{fig:PercentageUsedSimulationVSObserved}
\end{figure}

An important operational indicator is the percentage of utilization. 
Fig.~\ref{fig:PercentageUsedSimulationVSObserved} shows the utilization for each process, where (a) shows the values obtained from the simulation with the proposed model and (b) concerns to the historical orders. 
More than the 100 percent means that the extended capacity of certain process was used.
We observe that the utilization of resources is more regular in the model's simulation with few extra-times in which extended capacity is needed. 
Aside, the blanking process is the most used resource, meanwhile the annealing process is the loose resource.


\begin{figure}[h!]
\includegraphics[width=10cm]{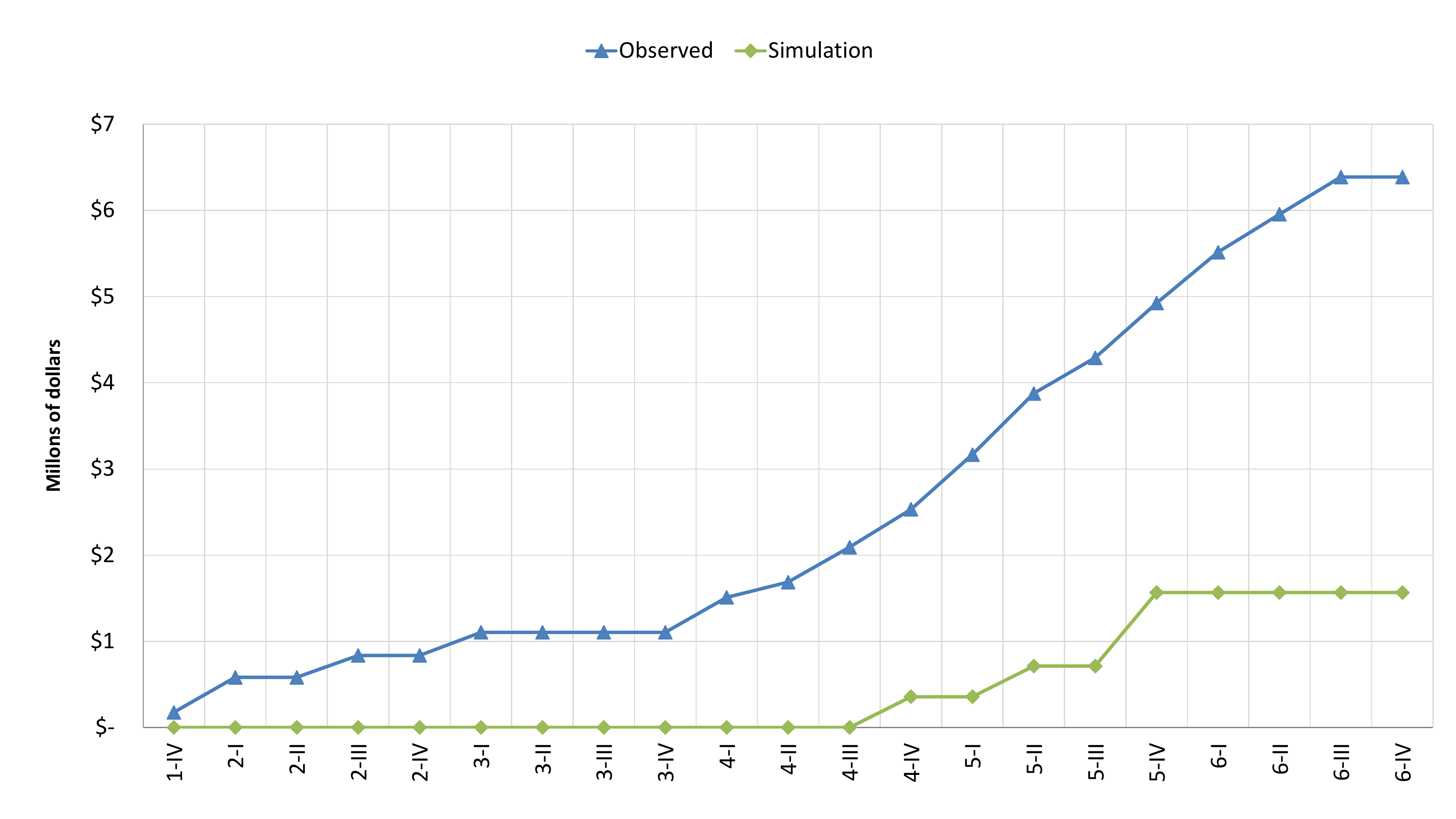}
\caption{Accumulated costs of the observed and simulation orders.}
\label{fig:CostSimulationVSObserved}
\end{figure}

We finish this section by showing the differences in the total cost between the observed orders and the simulation orders.
In Fig.~\ref{fig:CostSimulationVSObserved} the accumulated production costs of the extra-shifts are shown, the fixed costs are omitted.
The cost of the proposed model reduces in 24\% the costs of the historical minting orders.

\section{Conclusions}\label{sec_conclusion}
The issue of production cost of cash has become an important topic for central banks.
The central banks have developed several strategies to reduce the cost, they vary from the applying new technology in the fabrication, such as polymer, to reduce or enlarge the spacing between denominations.
This paper proposed a mixed-integer programming model to minimize the production costs of minting, while satisfying the production constrains, coin demand, and the minimum level of safety stocks.
Simulations were performed with data of 21 quarters yielding that this model could reduce the cost of the extra-shifts to 24\% compared to historical costs during 21 consecutive quarters. 

\section{Acknowledgement}
The first and the second authors were partially supported by SNI and CONACyT.

\end{document}